\documentclass{amsart}
\usepackage[all]{xy}
\usepackage{graphicx}
\usepackage{psfrag}
\usepackage{amsthm}
\usepackage{amscd}
\usepackage{amsfonts}
\usepackage{amstext}
\usepackage{amssymb}
\usepackage{amsmath}
\usepackage{amsxtra}
\usepackage{plain}
\usepackage[all]{xy}

\newtheorem{prop}{\bf Proposition}[section]
\newtheorem{cor}[prop]{{\bf Corollary}}
\newtheorem{lem}[prop]{{\bf Lemma}}
\newtheorem{thm}[prop]{{\bf Theorem}}

\numberwithin{equation}{section}

\newenvironment{rem}{{\bf Remark }}{ }
\newenvironment{defn}{{\bf Definition }}{ }

\newenvironment{pf}{{\bf proof }}{\qed\endtrivlist}

\newcommand{\Z}{\mathbb{Z} }
\newcommand{\N}{\mathbb{N} }
\newcommand{\C}{\mathbb{C} }
\newcommand{\R}{\mathbb{R} }

\newcommand{\ind}{\operatorname{ind}}
\newcommand{\Id}{\operatorname{Id}}
\newcommand{\A}{\mathcal{A}}

\newcommand{\Ch}{\operatorname{Ch}}

\newcommand{\End}{\operatorname{End}}

\newcommand{\id}{\operatorname{Id}}

\newcommand{\str}{\operatorname{str }}

\newcommand{\oddind}{\operatorname{odd-ind}}
\newcommand{\im}{\operatorname{Im}}
\newcommand{\Dom}{\operatorname{Dom}}

\long\def\symbolfootnote[#1]#2{\begingroup%
\def\thefootnote{\fnsymbol{footnote}}\footnote[#1]{#2}\endgroup}

\begin{document}

\title{Index theory and partitioning by enlargeable hypersurfaces}
\author{Mostafa ESFAHANI ZADEH}

\symbolfootnote[0]{\emph{2000 Mathematics Subject Classification}. 58J22 
(19K56 46L80 53C21 53C27).} 
\symbolfootnote[0]{\textit{Key words and phrases}. Higer index theory, 
enlargeablity, Dirac operators.}
\email{zadeh@uni-math.gwdg.de}

\address{Mostafa Esfahani Zadeh.\newline
 Mathematisches Institute, 
Georg-August-Universit\"{a}t,\newline
G\"{o}ttingen, Germany \and \newline 
Institute for Advanced Studies in Basic Siences(IASBS),\newline 
Zanjan-Iran}

\maketitle
\begin{abstract}
In this paper we state and prove a higher index theorem for 
an odd-dimensional connected 
spin riemannian manifold $(M,g)$ which is partitioned by an 
oriented closed hypersurface 
$N$. This index theorem generalizes a theorem due to N. Higson 
in the context of Hilbert 
modules. Then we apply this theorem to prove that if $N$ is area-enlargeable and if 
there is a smooth map from $M$ into $N$ such that its restriction to $N$ has 
non-zero degree then the the scalar curvature of $g$ cannot be uniformly positive.
\end{abstract}

\section{Introduction}
Given a compact manifold $M$, it is certainly an interesting problem to 
decide wether
it carries a riemannian metric with 
{\it everywhere positive scalar curvature} or not. 
This problem is revealed 
to be also very difficult. For constructing a metric with positive 
scalar curvature the most 
powerful technique is the Gromov-Lawson and Schoen-Yau surgery theorem asserting 
that if $M$ has a metric with positive scalar curvature 
and if $M'$ is obtained from $M$ by performing surgeries in co-dimension 
grater than or equal to $3$, then $M'$ carries a metric with positive 
scalar curvature too (see \cite{GrLa1} and \cite{ScYau}). 
In the other direction, i.e. to find obstruction for the existence of 
such metrics, the Atiyah-Singer index theorem and all its variant 
come into the play through the Lichnerowicz formula 
(see for example \cite{LaMi}). 
Even, it has been believed that all obstruction for the existence of 
metrics with positive scalar curvature on a spin manifold $M$ 
can be encapsulated in a sophisticated index $\alpha_{max}^\R(M)$ which takes 
its value in $KO_n(C^*_{max,\R}(G))$, where $G$ is the fundamental group 
of $M$ (see e.g. \cite{J.Ros3}). This assertion is 
known as the Gromov-Lawson-Rosenberg conjecture and is shown, by 
T. Schick in \cite{Schick-operator}, to not be true in this general form. 
Nevertheless, this index might subsume all index theoretic obstruction for 
the existence of metrics with positive scalar curvature on spin manifolds. 
One obstruction for the existence of metrics with positive scalar curvature is the 
enlargeability which was introduced by Gromov and Lawson 
(see \cite{GrLa2,GrLa3}). Enlargeability is a homotopy invariance 
of smooth manifolds and the category of enlargeable manifolds form a 
rich and interesting family containing, for example, all hyperbolic manifolds and 
all sufficiently large 3-dimensional manifolds. 

\begin{defn}\label{enlar}
Let $N$ be a closed oriented manifold of dimension $n$ with  
a fixed riemannian metric $g$. The manifold $N$ is enlargeable if for each 
real number $\epsilon>0$ there is a riemannian spin cover $(\tilde N,\tilde g)$, 
with lifted metric, and a smooth map $f:\tilde N\to S^n$ such that: 
the function $f$ is constant outside a compact subset $K$ of $\tilde N$; 
the degree of $f$ is non-zero;  and
the map $f:(\tilde N,~\tilde g)\to(S^n,g_0)$ is $\epsilon$-contracting, where
$g_0$  is the standard metric on $S^n$. Being $\epsilon$-contracting means that 
$\|T_xf\|\leq\epsilon$ for each $x\in \tilde N$, where 
$T_xf\colon T_x\tilde N\to T_{f(x)}S^n$. 
The manifold $N$ is said to be area-enlargeable if the function 
$f$ is $\epsilon$-area 
contracting. This means $\|\Lambda^2T_xf\|\leq\epsilon$ for each 
$x\in \tilde N$, where 
$\Lambda^2T_xf\colon \Lambda^2T_x\tilde N\to \Lambda^2T_{f(x)}S^n$. 
\end{defn}

It turns out that a closed area-enlargeable manifold cannot 
carry a positive scalar curvature and the basic tool to prove this theorem is a 
relative version of the Atiyah-Singer index theorem, c.f. 
\cite[theorem 4.18]{GrLa3}. So one may expect that the enlargeability 
obstruction be recovered by the index theoretic obstruction $\alpha_{max}^\R$. 
In fact T. Schick and B. Hanke in \cite{HaSc1,HaSc2} have 
proved that $\alpha_{max}^\R(N)\neq0$ if $N$ is enlargeable.

Given a complete riemannian manifold $(M,g)$ it is interesting to decide whether 
the scalar curvature of $g$ is uniformly positive. Beside its interest on 
itself, this question has clearly applications to the compact case too. 
The following result is the main result of 
this paper(see theorem \ref{thmpartenlar} is section \ref{sectionthree})  

{\it Let $(M,g)$ be a  complete riemannian manifold which has a finite spin covering,
and let $N$ be a 
closed area-enlargeable sub manifold of $M$ with co-dimension $1$. If there is a 
smooth map $\phi:M\to N$ such that its restriction to $N$ is of non-zero 
degree then the scalar curvature of $g$ cannot be uniformly positive.}

To prove this assertion 
we have put together some basic results and methods introduced by N.Higson, J.Roe,  
B.Hanke and T.Schick concerning index theory in the context of operator algebras. 
This result seems not easy to be obtained by means of the geometric methods 
of \cite{GrLa3}. 
So it shows also the efficiency of the operator algebraic index theory to prove 
results on the non-existence of metric with positive scalar curvature.

With the above notation, let $E$ be a Clifford bundle over $M$ and put $H=L^2(M,E)$. 
This a Hilbert space which is assumed be   
acted on by a Dirac type operator $D$. The operator $U=(D+i)(D-i)^{-1}$ is 
bounded on $H$. 
Let $N$ be a closed oriented hypersurface which partitions $M$ into two sub manifold 
$M_-$ and $M_+$ with common boundary $N$. 
The restriction of $D$ to $N$ defines a Dirac type operator $D_N$ 
with Fredholm index $\ind D_N$. Let $\phi_+$ denote a smooth function on $M$ 
which coincides to the characteristic function of $M_+$ outside a compact set and put 
$\phi_-=1-\phi_+$.
It turns out (see \cite{Higson-note}) that the bounded operator $U_+=\phi_-+\phi_+U$ 
is Fredholm and its index is denoted by $\ind(D,N)$. 
After an appropriate choice of orientations the following relation 
holds between this index and the Fredholm index of $D_N$, c.f. 
\cite[theorem 1.5]{Higson-note} and \cite[theorem 3.3]{Roe-partitioning} ,
\begin{equation}\label{higeq}
\ind(D,N)=\ind D_N~.
\end{equation}

This formula has two immediate applications. First one is to provide a 
proof for the cobordism invariant of the analytical index of a Dirac type operator, 
and the second one is the following: If $\hat A(M)\neq0$ then the scalar curvature 
of $g$ cannot be uniformly positive. We generalize the above theorem 
in the context of Hilbert Module over $C^*$-algebra, i.e. instead of $E$ we consider 
a Clifford Hilbert $A$-module bundle $W$. For example $W=S(M)\otimes V$ 
where $S(M)$ is the 
spin bundle of $M$ and $V$ is a Hilbert $A$-module bundle over $M$. In this case the 
twisted Dirac operator is denoted by  $D^V$ and its restriction to $N$ by $D_N^V$. The 
indices $\ind(D^V,N)$ and $D_N^V$ are elements of $K$-group $K_0(A)$ and 
we show in theorem 
\ref{higpar} of the section \ref{sectiontwo} that
\[\ind(D^V,N)=\ind D_N^V~.\]
As in above this relation can be used to prove the cobordism invariance 
of $\ind D_N^V$, 
c.f. corollary \ref{cobinv}. We have already mentioned the 
theorem \ref{thmpartenlar} in section 
\ref{sectionthree}. This theorem 
should be considered as a counterpart of the second application of 
the relation \eqref{higeq}. 

{\bf Acknowledgment :} The author would like to thank Thomas Schick for helpful 
discussions and John Roe for bringing the paper \cite{Higson-note} to his attention .

\section{Index theory on odd dimensional partitioned 
complete manifolds}\label{sectiontwo}

Let $(M,g)$ be an oriented complete non-compact manifold and let $W$ be a 
Clifford bundle on $M$ which is at the same time a Hilbert $A$-module bundle. 
We assume that this bundle is equipped with a connection which is compatible 
to the Clifford action of $TM$ and denote the corresponding $A$-linear Dirac 
type operator by $D$.  
For $\sigma$ and $\eta$ two compactly supported  smooth sections of $W$ put 
\[\langle\sigma,\eta\rangle=\int_M\langle\sigma(x),\eta(x)\rangle\,d\mu_g(x)\in A.\] 
It is easy to show that 
$|\sigma|=\|\langle\sigma,\sigma\rangle\|^{1/2}$ is a norm on $C_c(M,W)$. 
The completion of $C_c^\infty(M,W)$ with respect to this norm is the 
Hilbert $A$-module 
$H=L^2(M,W)$. The operator $D$ is a formally self adjoint operatot 
on this Hilbert module.

Let $T$ be an $A$-linear map which is defined on 
a dense subspace $\Dom(T)$ of $H$. 
The graph of $T$ is the following subset of $H\oplus H$
\[graph (T):=\{(u,T(u))|u\in \Dom(T)\}\]
The closure of this graph with respect to norm topology in 
$H\oplus H$ turns out to be 
the graph of an operator $\bar T$ which is called {the closure} of $T$. It is clear that 
$\Dom(\bar T)$ is a closed $A$-subspace of $H$.
The adjoint of $T$ is the closed operator 
$T^*$ such that $\langle Tu,v\rangle=\langle u,T^*v\rangle$ for $u\in\Dom(T)$ and $v\in\Dom(T^*)$. 
$T$ is said to be self-adjoint if $\Dom(T)=\Dom(T^*)$ and $T^*=T$. 
The operator $T$ is called regular if there is a bounded adjointable operator 
$P\in\mathcal L_A(H\oplus H,H)$ with $\im P=\Dom\bar T$ 
(c.f proposition 5 of \cite{Dan}).
If $T$ is regular then for each positive real number $c$ the image of the operator 
$c+T^*T$ is a dense sub space of $H$, c.f. \cite[proposition 4]{Dan}. 
In particular if $T$ is self adjoint then both 
$T+i$ and $T-i$ have dense images because $(T\pm i)(T\mp i)=T^2+1$ has dense image.
%Regularity and self adjointness of the operator $T$ make it possible to associate to 
%each continuous (not necessarily bounded) function $f$ on $\spec(T)$ a closed $A$-linear 
%operator $f(T)$ on $H$, c.f. \cite[proposition 16]{Dan}. 
%This correspondence define a continuous functional calculus for $T$. 
As a consequence of above discussion, by the following lemma, $(D\pm i)^{-1}$ is a 
a well defined bounded operator on $H$ 
 \begin{lem}\label{regularity}
 The operator $D$ is regular. 
\end{lem}
\begin{pf}
The operator $D$ is a formally self adjoint operatot 
on $H$, so following the above arguments it suffices to prove that the domain 
of $D$ is dense and that  $\Dom(\bar D))$ 
is the image of a bounded operator from $\mathbb H\oplus \mathbb H$ to $\mathbb H$. 
We recall that in \cite{FoMi} an $A$-linear pseudodifferential 
calculus is developed which provides a parametrix $Q$ for $D$, i.e. a 
pseudodifferential operator of order $-1$ so that 
\[DQ=\id+R~~\text{ and }~~QD=\id+S,\]
where $R$ and $S$ are smoothing operators and $Q$ and  
$S$ are bounded adjointable operators on $H$. Here the Dirac operator $D$ has   
domain $H^1(M,W)$, the $1$-th order Sobolev space which is a dense subspace of 
$H$. If $\eta=\bar D(\xi)$ then there is a sequence $\{\xi_k\}_k$ in $\Dom(D)$ such 
that $\eta=\lim_{k\to\infty} D(\xi_k)$ and $\lim_{k\to\infty}\xi_k\to\xi$. Then, since $Q$ and $S$ 
are bounded operator from $H$ to $H^1(M,W)$ we conclude  
$Q(\eta)=\xi+S(\xi)$, so $\xi=Q(\eta)-S(\xi)\in\im(Q)+\im(S)$, 
i.e. $Dom(\bar D)\subset \im Q+\im S$. 
Conversely, $Q$ and $S$ are pseudodifferential operators of order, respectively, $-1$ 
and $-\infty$, so they map $\mathbb H$ into $H^1(M,W)$. 
Therefore $Dom(\bar D)=\im (Q\oplus S)$ and this completes the proof.
\end{pf}

For bounded operators $P$ and $Q$ in $\End_A(H)$ 
by  $P\sim Q$ we mean that the difference $P-Q$ is compact in 
the sense of \cite{FoMi}. 
The following simple lemma is a key tool in what follows
\begin{lem}\label{keylem}
\begin{enumerate}
\item The operator $(D+i)^{-1}$ is a bounded operator on $H$ 
and $\theta(D+i)^{-1}\sim0$ 
where $\theta$ is a compactly supported function on $M$. 
\item If $\phi$ is a smooth function of $M$ which is locally constant 
outside a 
compact subset then \[[(D\pm i)^{-1},\phi]\sim0~.\]
\end{enumerate}
\end{lem}
\begin{pf}
For $\sigma$ a smooth element of $H$ we have 
\begin{align*}
\langle (D+i)\sigma,(D+i)\sigma\rangle&=\langle(D^2+1)\sigma,\sigma\rangle\\
&\geq\langle\sigma,\sigma\rangle~,
\end{align*}
therefore $\|(D+i)^{-1}|\leq1$, i.e. 
$(D+i)^{-1}$ is bounded. 

The operator $(D+i)$ is an elliptic differential $A$-operator in the sense 
of \cite{FoMi}, 
so there is pseudodifferential $A$-operators $Q$ and $S$ of negative order 
such that 
\[Q(D+i)-\id=S~,\]
which gives rise to the following relation 
\[\theta Q-\theta(D+i)^{-1}=\theta S(D+i)^{-1}~.\] 
This shows that it is enough to prove the compactness of the 
operators $\theta Q$ and $\theta S$. 
Moreover for each positive number $\delta$ we may assume  $Q$ 
and $S$ be of $\delta$-propagation speed. This can be done by multiplying 
the kernels of this operator with a bombe function which is supported 
around the diagonal of $M\times M$.
Given a bounded sequence $\{\gamma_j\}$ in $H$ we show that the sequence 
$(Q\gamma_j)_{|supp\,\theta}$ contains a convergence sub-sequent and 
similarly for the sequence $(S\gamma_j)_{|supp\,\theta}$. For this purpose 
let $M'$ be an open sub manifold of $M$ with compact closure such that it 
contains $2\delta$-neighborhood of $supp(\theta)$. Deform all geometric structure 
in the 
$\delta$ neighborhood of $\partial M'$ to product structures and consider the 
double compact manifold  $(M'\sqcup_{\partial M'}M'^-)$ with Clifford Hilbert 
$A$-module bundle induced from $W$. 
Coresponding to the sequence $\{\gamma_j\}$ in $H$ consider the sequence 
$\{\gamma'_j\}$ in 
$(M'\sqcup_{\partial M'}M'^-)$ such that $\gamma'_{j|M'}=\gamma_{j|M'}$. 
Since $Q$ has propagation speed $\delta$  
\begin{equation}\label{dabsr}
(Q\gamma_j)_{|supp\,\theta}=(Q'\gamma'_j)_{|supp\,\theta}~.
\end{equation}
Here $Q'$ is any pseudodifferential $A$-operator of negative order 
on double space whose kernel is 
equal to the kernel of $Q$ on $M'$. Being of order negative, $Q'$ is a 
compact operator on $H$ (c.f. \cite[page 109]{FoMi}) so $\{Q'\gamma_j\}$ has 
a $L^2$-convergent subsequent. The restriction of this sub-sequent to $supp\,\theta$ 
is convergent too and this prove what we wanted in view of the relation \eqref{dabsr}.
The compactness of $\theta S$ can be proved completely similar. 

To prove the second part notice that 
\[[(D\pm i)^{-1},\phi]=(D\pm i)^{-1}[\phi,D](D\pm i)^{-1}~.\]
On the other hand $[D,\phi]\sigma=\rm{grad}(\phi).\sigma$ where "." denote 
the Clifford action. Since the gradient vector field $\rm{grad}(\phi)$ vanishes 
outside a compact set, for an appropriate compactly supported function $\theta$ 
on $M$ one has $[D,\phi]=[D,\phi]\theta$, so the first part of the lemma can be 
used to deduce the compactness of $[\phi,D](D\pm i)^{-1}$ and 
hencefore the compactness of $[(D\pm i)^{-1},\phi]$
\end{pf}

Suppose now that $M$ is partitioned by an oriented compact hypersurface $N$ into two 
parts $M_+$ and $M_-$. We assume that the the positive direction of the unite normal 
vector $\vec n$ to $N\subset M$ points out from $M_-$ to $M_+$. Let $\phi_+$ 
be a smooth function on $M$ 
which is equal to the characteristic function of $M_+$ outside a compact subset 
of $M$, and put $\phi_-=1-\phi_+$. To state the next theorem we recall 
from \cite[page 96]{FoMi} that a bounded operator $P\in End_A(H)$ is a 
Fredholm $A$-operator if there is a decomposition $H=H_0\oplus H_1$ of the source 
space and a decomposition $H=H'_0\oplus H'_1$ of the target space such that $H_0$ 
and $H'_0$ are finitely generated $A$-modules and if the matrix form of $P$, 
with respect to these decompositions, is given by 
\[\left(
\begin{array}{cc}
P_0&0\\0&P_1
\end{array}
\right)\]
where $P_1\colon H_1\to H'_1$ is isomorphism with bounded inverse. 
The Fomenko-Mischenko index of the Fredholm $A$-operator $P$ is given by 
\[\ind P=[H_0]-[H'_0]\in K_0(A)\]
and turns out to be independent of all choices made in its definition. This 
index is a homotopy invariant and is invariant with respect to perturbations 
by compact 
$A$- operators, c.f. lemmas $1.5$ and $2.3$ in \cite{FoMi}. Morover if $Q$ 
is another Fredholm $A$-operator then (c.f. \cite[lemma 2.3]{FoMi})
\[\ind(PQ)=\ind P+\ind Q\in K_0(A)~.\]
\begin{thm}\label{defther}
Let $U=(D-i)(D+i)^{-1}$ be the Cayley transform of $D$. We have 
\begin{enumerate}
\item the operators $U_+=\phi_-+\phi_+U$ and $U_-=\phi_++\phi_-U$ are 
Fredholm operators in $End_A(H)$
\item $\ind(U_+)=-\ind(U_-)$;
\item The value of $\ind(U_+)$ does not depend on $\phi_+$ but on the 
cobordism class of the partitioning manifold $N$. 
\end{enumerate}
\end{thm}
\begin{pf}
Clearly both $U_\pm$ are $A$-linear and bounded. Since $U=1-2i(D+i)^{-1}$ 
the second part of the above lemma implies 
\begin{align*}
\phi_+U&=\phi_+-2i\phi_+(D+i)^{-1}\\
&\sim \phi_+-2i(D+i)^{-1}\phi_+=U\phi_+
\end{align*}
and similarly $\phi_-U\sim U\phi_-$. Since the support of 
$\phi_+\phi_-$ is compact, it follows from the first part of the previous lemma that 
\[\phi_+\phi_-U=\phi_+\phi_--2i\phi_+\phi_-(D+i)^{-1}\sim\phi_+\phi_-~.\]
Using these relation one has $U_\pm U_\pm^*\sim\id$ and $U_\pm^* U_\pm\sim\id$, 
so both 
$U_+$ and $U_-$ are Fredholm $A$-operator according to \cite[theorem 2.4]{FoMi}. 
Since $U_+U_-=\id$ we get the following relation 
\[\ind U_++\ind U_-=\ind U=0\in K_0(A)~.\]
A different choice for $\phi_+$ is of the form $\phi_++\phi$ where 
$\phi$ is a compactly supported smooth function. The corresponding  operator 
differs from $U_+$ by 
\[-\phi+\phi U=-\phi+\phi(1-2i(D+i)^{-1})\sim0~,\]
where the equivalence is coming from first part of the previous lemma. 
Consequently the index of $U_+$ does not depend on the choice of $\phi_+$. 
The third part of the theorem is a direct consequence of the second part. 
\end{pf}

The Fomenko-Mischenko index of the  operator $U_+$ is denoted by 
$\ind(D,N)\in K_0(A)$. 
The following property of this index is crucial for our purposes

\begin{thm}\label{asliposi}
If $D$ is an isomorphism, i.e. has the bounded inverse, then $\ind(D,N)=0$~. 
In particular 
if $D^V$ is the spin Dirac operator on $M$ twisted by the 
Hilbert $A$-module bundle $V$, if the scalar curvature $\kappa$ of $g$ is 
uniformly positive 
and the curvature of $V$ is sufficiently 
small then $\ind(D^V,N)=0\in K_0(A)$. 
\end{thm}
\begin{pf}
If $D$ has bounded inverse then for $0\leq s\leq1$ the families 
$(D-si)$ and $(D+si)^{-1}$ of bounded operator are continuous. 
So $U(s):=(D-si)(D+si)^{-1}$ is a homotopy of bounded operators 
and $U_+(s)=\phi_-+\phi_+U(s)$ is 
a homotopy of Fredholm $A$-operators between $U_+$ and $\id$. 
By the homotopy invariance of the the Fomenko-Mishenkov index we get $\ind(D,N)=0$. 

To prove the second part 
let $\kappa\geq\kappa_0>0$. By the Lichnerowicz formula  
\[(D^V)^2=\nabla^*\nabla+\frac{1}{4}\kappa+ R~,\]
where $ R\in \End_A(V)$ is sufficiently small, 
say $\| R\|^\leq \frac{1}{8}\kappa_0$.
So we get the following inequality in $A$, where $\sigma$ is a smooth 
$L^2$-section of $S(M)\otimes V$
\[\langle D^V(\sigma),D^V(\sigma)\rangle\geq
\frac{1}{8}\kappa_0\langle\sigma,\sigma\rangle\]
which implies the boundedness of $(D^V)^{-1}$. 
The assertion follows now from the first part. 
\end{pf}
\begin{rem}
The proof of the above theorem may be slightly modified to show that $\ind(D,N)=0$ if 
there is a gap in the $L^2$-spectrum of $D$. 
\end{rem}

The following lemma shows that the index $\ind(D,N)$ is invariant with respect to 
modifications of data at each partition $M_+$ or $M_-$.

\begin{lem}\label{nimiso}
For $j=1,2$ let $M_j$ be a complete manifold partitioned by compact hypersurface 
$N_j\subset M_j$ and let $W_j$ be a Clifford Hilbert $A$-module bundle on $M_j$. 
If there is an isometry  $\gamma\colon M_{2+}\to M_{1+}$  
which is lifted to an isomorphism of 
Clifford and Hilbert module structures then 
\[\ind(D_1,N_1)=\ind(D_2,N_2).\]
The similar assertion is true for $M_{j-}$.
\end{lem}
\begin{pf}
Let $\phi_1$ be a smooth function on $M_1$ which vanishes in a neighborhood of 
$M_{1-}$ and is equal to $1$ outsid a compact subset of $M_{1+}$. Notice 
that  $\phi_2=\phi_1\circ\gamma$ is defined only on $M_{2+}$ but it can 
be extended by zero 
to whole $M_2$. As in the theorem  \ref{defther} we have 
\begin{gather*}
U_{1+}=1+2i\phi_1(D_1+i)^{-1}\\
U_{2+}=1+2i\phi_2(D_2+i)^{-1}\sim 1+2i(D_2+i)^{-1}\phi_2~.
\end{gather*}
The isomorphism $\gamma$ provides an unitary isomorphism 
$\Gamma\colon L^2(M_{1+},W)\to L^2(M_{2+},W)$. By taking an arbitrary isomorphism 
$L^2(M_{1-},W)\to L^2(M_{2-},W)$ and using the direct decomposition 
$L^2(M_j,W)=L^2(M_{j-},W)\oplus L^2(M_{j+},W)$ we get an isomorphism 
\[T:L^2(M_1,W)\to L^2(M_2,W)~.\]
One has $(D_2+i)\Gamma\phi_1=\Gamma(D_1+i)\phi_1$ 
and $\phi_2\Gamma(D_1+i)=\Gamma\phi_1(D_1+i)$, so
\begin{align*}
TU_{1+}-U_{2+}T&\sim T(1+2i\phi_1(D_1+i)^{-1})-(1+2i(D_2+i)^{-1}\phi_2)T\\
&=2i(\Gamma\phi_1(D_1+i)^{-1}-(D_2+i)^{-1}\phi_2\Gamma)\\
&=2i(D_2+i)^{-1}\Gamma[D_1,\phi_1](D_1+i)^{-1}
\end{align*}
Now one can proceed as in the proof of the second part of the lemma \ref{keylem} to 
deduce that the last expression is a compact operator. Consequently 
$TU_{1+}T^{-1}\sim U_{2+}$ which implies  
\[\ind U_{1+}=\ind U_{2+}\in K_0(A)~.\] 
\end{pf}

The Clifford action of $i\vec n$ provides a $\Z_2$ grading for $W_{|N}$. 
Let $D_N$ denote the Dirac type operator acting on smooth sections of $W_{|N}\to N$. 
This is a $A$-linear elliptic operator and has an index $\ind D_N\in K_0(A)$. 

The following theorem, as well as its proof, is a generalization of 
the theorem $1.5$ of 
\cite{Higson-note} to the context of the Hilbert module bundles and 
its proof follows the same 
(see also \cite[]{Roe-partitioning}). 
\begin{thm}\label{higpar} The following equality holds in the $K$-group $K_0(A)$
\begin{equation*}
\ind D_N=\ind(D,N)~.
\end{equation*}
\end{thm}
\begin{pf}
As the first step we show that it is enough to prove the theorem for the cylindrical 
manifold $\R\times N$ with product metric $(dx)^2+g_N$ and pull back bundle 
$p^*(V_{|N})$ where $p$ is the projection of $\R\times N$ onto the second factor. 
Consider a collar neighberhood $(-1,1)\times N$ in $M$. 
Using the lemma \ref{nimiso} we may change $M_-$ to the product 
form $(-\infty,1/2)\times N$ without 
changing the index $\ind(M,N)$. 
By applying the third part of the theorem \ref{defther} we may assume 
that the partitioning manifold 
is $\{1/2\}\times N$, then the lemma \ref{nimiso} can be used again to 
replace $M$ with the 
cylinder $\R\times N$ without changing $\ind(M,N)$ nor $\ind D_N$. Consequently, 
to prove the theorem it suffices to prove it in the special case of the cylinder. 
At first we prove the theorem for the very special case of the  
euclidian Dirac operator $-id/dx$ twisted by the finitely projective $A$-module $V_0$. 
We denote this twisted Dirac operator by $D_\R^{V_0}$ and prove the following 
relation  
\begin{equation}\label{dfrash}
 \ind(\pm D_\R^{V_0})=\pm [V_0]\in K_0(A)~.
\end{equation}
For this purpose let $\phi_+$ be a smooth function on $\R$ satisfying the conditions of 
the theorem \ref{defther} and put $\psi=2\phi_+-1$. One has
\begin{equation*}
 U_+=(D_\R^{V_0}-i\psi)(D+i)^{-1}~.
\end{equation*}
Therefore the $L^2$-kernels of $U_+$ and $U_+^*$, as $A$-modules, are isomorphic to $\ker(D-i\psi)$ 
and $\ker(D+i\psi)$. The space of $L^2$ solutions of $U_+^*=-id/dx+i\psi$ is null while 
the space of the $L^2$-solution of  $U_+=-id/dx-i\psi$ consists of the following 
smooth functions 
\[f(x)=\exp(-\int_0^x\psi(t)\,dt)v~;\hspace{1cm}v\in V_0.\] 
Consequently the $L^2$-kernels of $U_+$ and $U_+^*$ are isomorphic to the finitely 
generated projective $A$-modles $V_0$ and $0$, so $\ind(U_+)=[V_0]-0$ which is the 
desired relation. The case of $-D_\R^{V_0}$ is similar. 

Now we are going to prove the theorem for the 
cylinder $(\R\times N, (dx)^2+g_N)$, which completes the proof as is explained in above. 
The operator  $D^V$ has the following form 

\begin{equation*}
D^V=\left( 
\begin{array}{cc}
i\partial_x&(D_N^V)^-\\
(D_N^V)^+&-i\partial_x
\end{array}\right). 
\end{equation*}
Consider the Dirac type operator $D_N^V$
as an unbounded operator on $L^2(N, W_{|N})$. As it has been pointed out 
just before the theorem \ref{defther}, there is a 
decomposition of $L^2(N,W_{|N})=\mathcal W_0\oplus\mathcal W_1$ to direct 
sum of invariant $A$-modules such that $\mathcal W_0$ is finitely  
generated and projective. 
Moreover the restriction of $D_N^W$ to $\mathcal W_1$ has a bounded inverse. 
The operator $\mathcal R:=(D_N^V)_{|\mathcal W_0}\oplus 0$ is a compact 
$A$-operator on $L^2(N,W_{|N})$, so the following operator 
\[\tilde D_N:=\left( 
\begin{array}{cc}
0&0\\0&(D_N^V)_{|\mathcal W_1}
\end{array}
\right)\]
is a compact perturbation of $D_N^V$, consequently 

\begin{equation}\label{khash0}
\ind D_N^V=\ind \tilde D_N=[\mathcal W_0^+]-[\mathcal W_0^-]\in K_0(A)~.
\end{equation} 
The Cayley transform of the following family of operators 
\[D_s=\epsilon i\frac{d}{dx}+
\left( 
\begin{array}{cc}
s\mathcal R&0\\0&(D_N^V)_{|\mathcal W_1}
\end{array}
\right)\]
is a continuous family of bounded operators to which the lemma 
\ref{keylem} is applicable. Here $\epsilon$ stands for the grading operator.
Therefore we have the homotopy $U_+(D_s)$ 
of Fredholm operators with the same index in $K_0(A)$. Therefore 

\begin{equation}\label{khash1}
\ind(D^V,N)=\ind(\bar D,N)~,
\end{equation}
where 
\[\bar D=\epsilon i\frac{d}{dx}+\left( 
\begin{array}{cc}
0&0\\0&(D_N^V)_{|\mathcal W_1}
\end{array}
\right)~.\]
\begin{rem}
Notice that $\bar D$ is not a twisted Dirac operator acting on sections of a 
Hilbert module bundle. Nevertheless since the lemma \ref{keylem} is applicable to this operator 
$U_+(\bar D)$ is a Fredholm $A$-operator on $L^2(M,\mathcal W_1)$ and we can define the index $\ind(\bar D,N)=\ind U_+(\bar D)$ as an element in $K_0(A)$.
\end{rem}
%By the above discussion the following relation will prove the theorem 
%\ref{higpar} 
%in the special case of the cylinder.  
%\begin{equation}\label{amimo}
%\ind(\bar D,N)=[\mathcal W_0^+]-[\mathcal W_0^-]\in K_0(A)~,
%\end{equation}
%To prove this relation
With respect to  direct sum 
$L^2(\R\times N,W)=L^2(\R,\mathcal W_0)\oplus L^2(\R,\mathcal W_1)$ we put 
$\bar D=\bar D_0\oplus \bar D_1$. In the view of the above remark it is clear that  
\[\ind(\bar D,N)=\ind(\bar D_0,N)+\ind(\bar D_1,N)\in K_0(A).\] 
Notice that $\bar D_0$ is the euclidian Dirac operator (up to sign) twisted by a 
finitely generated module, so  $\ind(\bar D_0,N)$ is well defined. 
If $\sigma$ is a smooth element of $L^2(\R,\mathcal W_1)$ then 
\begin{align*}
\|(\bar D_1\pm i\psi)\sigma\|^2&=\langle(\bar D_1\mp i\psi)(\bar D_1\pm i\psi)\sigma,\sigma)\rangle\\ 
&=\langle(D_N^V)^*D_N^V\sigma+((id/dx\pm i\psi)^*(id/dx\pm i\psi))\sigma)\rangle\\
&\geq\|D_N^V\sigma\|^2\\
&\geq\delta\|\sigma\|^2~,\hspace{5mm}\text{ for a } \delta>0
\end{align*}
where the last inequality results from the fact that $(D_N^V)_{|\mathcal W_1}$ 
has a continuous inverse. Now the argument in the proof of the theorem  can 
be applied to deduce $\ind(\bar D_1,N)=0\in K_0(A)$.  
On the other hand, the operator $\bar D_0$ has the following form acting on 
$L^2(\R,\mathcal W_0^+)\oplus L^2(\R,\mathcal W_0^-)$
\[D_\R^{\mathcal W_0}\oplus-D_\R^{\mathcal W_1}\]
Therefore the following equality follows by applying the relations \eqref{dfrash}  
\begin{equation*}
\ind(\bar D_0,N)=[\mathcal W_0^+]-[\mathcal W_0^-]\in K_0(A)~.
\end{equation*}
This relation along relations \eqref{khash0} and \eqref{khash1} prove the 
theorem for the cylindrical case and so complete the proof of the theorem.
\end{pf}
The following theorem is an immediate application of the previous theorem. 
\begin{cor}\label{cobinv}
Let $N$ be a closed even dimensional manifold and $W$ be a Clifford Hilbert 
$A$-module bundle on $N$. Let $D$ be a Dirac type operator acting on sections of $W$. 
If there is a compact manifold $M$ with $N=\partial M$ and if all geometric 
structures extend to $M$, then $\ind D=0\in K_0(A)$
\end{cor}

\section{Partitioning by enlargeable manifolds}\label{sectionthree}

In this section we apply theorem \ref{higpar} to prove the following theorem 
concerning the existence of complete metrics on non-compact manifolds 
with uniformly positive scalar curvature.

\begin{thm}\label{thmpartenlar}
Let $(M,g)$ be a non-compact orientable complete $n$-manifold where 
$n\geq2$ and let $M$
has a finite covering which is spin. 
Let $N$ be a $(n-1)$-dimensional area-enlargeable closed sub-manifold of $M$ which 
 partitions $M$. If there is a map $\phi\colon M\to
N$ such that its restriction to $N$ has non-zero degree, then the
scalar curvature of $g$ cannot be uniformly positive. 
\end{thm}

\begin{pf}
At first notice that $M$, hence $N$, may be assumed to be spin. If not we 
consider the 
finite spin riemannian covering $(\tilde M,\tilde g)$ where $\tilde g$ is 
the lifting of $g$. Let $\tilde N\subset \tilde M$ be the induced covering for 
$N$ which is area-enlargeable. Its normal bundle is trivial so 
$\tilde N$ is spin too. 
The function $\phi$ has a lifting to a function 
$\tilde\phi:\tilde M\to\tilde N$ such that its restriction to 
$\tilde N$ is of non-zero degree. Moreover, the scalar curvature of $g$ is 
uniformly positive if and only if the 
scalar curvature of $\tilde g$ is uniformly positive. 
This discussion show that if we prove the assertion of the theorem for 
$\tilde M$, then it follows for $M$, too.

We can also assume $n$ to be an odd integer. If not consider 
the complete manifold $(M\times S^1,~g\oplus g_0)$ where $g_0$ is 
any riemannian metric on $S^1$. 
If $M$ has a finite spin covering $\tilde M$ then $\tilde M\times S^1$ 
is a finite spin covering for $M\times S^1$.
The restriction of the map 
$\phi\times \Id:M\times S^1\to N\times S^1$ has non-zero degree and 
$N\times S^1$ is area-enlargeable. Moreover if the scalar curvature of $g$ 
is uniformly positive then the scalar 
curvature of $g\oplus g_0$ is uniformely positive, so it suffices to prove the theorem for 
the odd dimensional complete manifold $(M\times S^1,~g\oplus g_0)$.

To prove this theorem we use some methods and constructions introduced by 
B.Hanke and T.Schick in \cite{HaSc1}  and \cite{HaSc2}. 
Following \cite[proposition 1.5]{HaSc2}, since $N$ is area-enlargeable, 
for each positive natural number $k$ there exists a $C^*$ algebra $A_k$ 
and a Hilbert $A_k$-module bundle $V_k\to N$ with connection $\nabla_k$ with 
following properties: The curvature $\Omega_k$ of $V_k$ satisfies

\[\|\Omega_k\|\leq\frac{1}{k},\]
and there exists a split extension
\[0\to\mathbb K\to A_k\to\Gamma_k\to0~,\]
where $\mathbb K$ denotes the algebra of compact operator on an infinite 
dimensional and separable Hilbert space, and $\Gamma_k$ is a certain 
$C^*$-algebra. In particular 
each $K_0(A_k)$ canonically splits off a $\Z=K_0(\mathbb K)$ summand. Denote by $W_k$
the spin bundle twisted by $V_k$ and by $D_N^{V_k}$ the associated twisted Dirac operator.
If $a_k\in K_0(A_k)$ denotes the index of $D_N^{V_k}$ then 
the $\Z=K_0(\mathbb K)$-component $z_k$ of $a_k$ 
is non-zero. 
Moreover there is a dense subalgebra 
$\A_k$ of $A_k$ which is closed under holomorphic calculus and there is a 
continuous trace $\alpha_k\colon \A_k\to\C$
such that $z_k=\ind_{\alpha_k} D_N^{V_k}$. Here we have used the fact that 
$\ind D_N^{V_k}\in K_0(A_k)=K_0(\mathcal A_k)$ so the expression 
$\ind_{\alpha_k} D_N^{V_k}$ makes sens. This index can be calculated in term of 
the geometry of $(M,g)$ and of the bundle $(V_k,\nabla^{V_k})$.
The theorem 9.2 of \cite{Sc-L^2} 
gives the following explicite formula for this index

\begin{equation}\label{expind}
z_k=\ind_{\alpha_k} D_N^{V_k}=\int_N A(TM)\wedge[\Ch_{\alpha_k}(V_k,\nabla^{V_k})]_+~.
\end{equation}
Here $[\omega]$ denotes the positive degree part of differential form 
$\omega\in\Omega^*(N)$ and  $\Ch_{\alpha_k}(V_k,\nabla^{V_k})$ is defined in term of the curvature 
$\Omega_k$ by the following relation which provides a closed 
differential form on $N$ 

\begin{equation}\label{relsum}
\Ch_{\alpha_k}(V_k,\nabla^{V_k})
:=\alpha_k(\str\sum_{j=0}^\infty\dfrac{\Omega_k\wedge
\dots\wedge\Omega_k}{j!} )~.
\end{equation}
The class of $\Ch_{\alpha_k}(V_k,\nabla^{V_k})$ in the de-Rham cohomology 
of $N$ is independent of the 
connection $\nabla^{V_k}$, so it is determined by the class of $V_k$ 
in $K_0(A_k)$. Since $K_0(\mathcal A_k)\simeq K_0(A_k)$ one can assume that 
the value of the expression between the parentheses at the right hand side of 
\eqref{relsum} belons to $\mathcal A^{ab}$, so in the domain of $\alpha_k$. 
This justifies the definition of the Chern character.

It is clear from the definition that for a smooth function $\psi\colon N\to N$ one 
has 
\[\Ch_{\alpha_k}(\psi^*V_k,\nabla^{\psi^*V_k})=\psi^*\Ch_{\alpha_k}(V_k,\nabla^{V_k})~.\]
The main feature of the virtual bundle $V_k$, coming from 
area-enlargeability of $N$, is that the Chern character 
$\Ch_\gamma(V_k,\nabla^{V_k})$ is, in fact, a closed differential $(n-1)$-form. 
 So given a smooth map $\psi\colon N\to N$, the explicite 
formula \eqref{expind} implies the following relation  

\begin{equation}\label{kateb}
\ind_{\alpha_k}D_N^{\psi^*V_k}=\deg(\psi).\ind_{\alpha_k} D_N^{V_k}~.
\end{equation}

For being able to use theorem \ref{higpar}, we need 
to work with flat bundles on $N$. 
For this purpose we use another fundamental construction introduced 
in \cite{HaSc1}. This construction consists of {\it assembling} the algebras $A_k$, the 
almost flat sequence of 
bundles $V_k$ and the connections $\nabla^{V_k}$ to 
construct a $C^*$ algebra $A$, a Hilbert $A$-module bundle $V$ and a flat 
connection $\nabla^V$ such that the index of the twisted Dirac operator $D_N^V$ 
acting on smooth sections of $W=S\otimes V$, taking its value in $K_0(A)$,  
keeps the track of the index theoretic information of $D_N^{V_k}$ when $k$ 
goes toward infinity. 
Denote by $\prod^b A_k$ the $C^*$-algebra consisting of all 
uniformly bounded sequences $(a_1,a_2,\dots)$ with $a_k\in A_k$ and by 
$\prod^0 A_k$ the $C^*$-algebra consisting of all 
sequences $(a_1,a_2,\dots)$ such that the sequence $\{\|a_k\|\}_k$ converges to $0$. 
The $C^*$-algebra $A$ is defined by the following quotient 
\[A:=\dfrac{\prod^b A_k}{\prod^0 A_k}.\]
The $C^*$-algebras $\Gamma$ and $\mathcal K$ are constructed from $\{\Gamma_k\}_k$ and 
from $\{\mathbb K\}_k$ by semilar quotients. 
Clearly one has the following split exact sequence 
\begin{equation}\label{dans}
0\to\mathcal K\to A\to\Gamma\to0,
\end{equation}
which gives rise to the following commutative diagram with spliting rows 

\begin{equation}\label{eq:cstardiagram}
0\to K_0(\mathcal K)\to K_0(A)
\to K_0(\Gamma)\to0 
\end{equation}
It turns out that  
\[K_0(\mathcal K)\simeq
\dfrac{\prod K_0(\mathbb K)}{\bigoplus K_0(\mathbb K)}\simeq
\dfrac{\prod\Z}{\bigoplus\Z}.\]
The proposition 1.5 of \cite{HaSc1} and the discussion after that show that 
the component of $\ind D_N^V$ is $K_0(\mathcal K)$ can be represented, with 
respect 
to the above isomorphism, by 
$z=[(z_1,z_2,\dots)]$ where $z_k=\alpha_k(\ind D_N^{V_k})\neq0$ for all $k\in \N$. 
This implies 
the nonvanishing result $\ind D_N^V\neq0\in K_0(\mathcal K)$. 

For $\psi$ a smooth function on $N$ it turns out from the construction of the 
bundle $V$ that the pull back bundle $\psi^*(V)$ may be constructed 
by assembling the bundles $\psi^*(V_k)$. Using  \eqref{kateb} and the above 
description of the  $K_0(\mathcal K)$-component of the higher index 
$\ind(D_N^{\psi^*V})\in K_0(A)$ we conclude that the $K_0(\mathcal K)$-component of 
$\ind D_N^{\psi^*V}$ equals to $\deg(\psi)$-times of the 
$K_0(\mathcal K)$-component of $\ind D_N^V$. 
Since this last component in non-vanishing, we conclude 
In particular 
\begin{equation}\label{fanz}
\ind D_N^{\psi^*V}\neq0\in K_0(\mathcal A)~,
\end{equation}
provided that $\deg(\psi)\neq0$. 
Now we are able to apply theorem \ref{higpar}. Using the map $\phi:M\to N$ we 
construct the pull back bundle $\phi^*V$ and the pull-back 
connection  $\phi^*\nabla$ which is flat. Let $D^{\phi^*V}$ be 
the spin Dirac operator of  $M$ twisted by the flat Hilbert 
$A$-module bundle$(\phi^*V,\phi^*\nabla)$.  
The restriction of this bundle to $N$ is $\psi^* V$, where $\psi:=\phi_{|N}$ 
is of non-zero degree. By theorem \ref{higpar} we have   
\begin{equation*}
\ind(D^{\phi^*V},N)=\ind D_N^{\psi^*V}~,
\end{equation*}
which gives the following non-vanishing result by relation \eqref{fanz}
\begin{equation}\label{satar}
\ind(D^{\phi^*V},N)\neq0\in K_0(A).
\end{equation}

On the other hand, by lemma \ref{asliposi} if the scalar curvature of $g$ 
is uniformly positive then  $\oddind D^{\phi^*V}$ vanishes. This is 
in contradition with the above non-vanishing result and completes 
the proof of the theorem.
\end{pf}

\end{document}